%% file: lsi-ggm.tex
\documentclass [a4paper,twoside,11pt]{article}

\usepackage[latin1]{inputenc}
\usepackage{amsmath,amsfonts,amssymb}
\usepackage{vmargin,graphicx,theorem}
\usepackage[english]{babel}
\usepackage{enumerate}

\setpapersize[portrait]{A4}
\setmarginsrb{2.2cm}{2cm}{2.2cm}{1cm}{0cm}{0.1cm}{0.5cm}{2cm}

\include{macros}

\newcommand{\affil}[1]{{\small\sl #1}}
\newcommand{\email}[1]{{\small E-mail: {\textsf {#1}}}}
\newcommand{\http}[1]{{\small Internet: {\textsf {#1}}}}

\begin{document}

\title{\sl Modified logarithmic Sobolev inequalities and transportation inequalities }
\author{
  Ivan Gentil, Arnaud Guillin\\
\affil{Ceremade (UMR CNRS no. 7534), Universit\'e Paris IX-Dauphine,}\\
\affil{Place de Lattre de Tassigny, 75775 Paris C\'edex~16, France}\\
\email{\{gentil,guillin\}@ceremade.dauphine.fr}\\
\http{http://www.ceremade.dauphine.fr/
\raisebox{-4pt}{$\!\!\widetilde{\phantom{x}}$}\{gentil,guillin\}/}\\
\\
Laurent Miclo\\
\affil{Laboratoire de Statistique et Probabilit\'es, (UMR CNRS no. 5583),} \\
\affil{Universit\'e Paul-Sabatier,} \\
\affil{118 route de Narbonne, 31062 Toulouse C\'edex 4, France}\\
\email{miclo@math.ups-tlse.fr}\\
}
\date{\today}\maketitle\thispagestyle{empty}

\begin{abstract}
We present a class of modified logarithmic Sobolev inequality, interpolating between
Poincar\'e and logarithmic Sobolev inequalities,
suitable for measures of the type $\exp(-|x|^\al)$ or more complex $\exp(-|x|^\al\log^\beta(2+|x|))$
($\al\in]1,2[$ and $\be\in\dR$)
which lead to new concentration
inequalities. These modified inequalities share common properties with usual logarithmic
Sobolev inequalities, as
tensorisation or perturbation, and imply as well Poincar\'e inequality. We also study the
link between these new
modified logarithmic Sobolev inequalities and transportation inequalities. 
\end{abstract}
%

\section{Introduction}
\label{sec-intro}

A probability measure $\mu$ on
 $\dR^n$ satisfies a logarithmic Sobolev inequality if there exists $C<\infty$ such that,  for every
 smooth enough functions $f$ on $\dR^n$,
\begin{equation}
\label{eq-gross}
\ent{\mu}{f^2} \leq  C \int \ABS{\nabla f}^2 d\mu,
\end{equation}
where
$$ \ent{\mu}{f^2} = \int f^2 \log f^2 d\mu- \int f^2 d\mu  \log \int  f^2 d\mu $$
and where $\ABS{\nabla f}$ is the Euclidean length of the gradient $\nabla f$
of $f$.

\medskip

Gross in~\cite{gross} defines this inequality and shows   that
the canonical Gaussian measure with density $(2\pi )^{-n/2} {e}^{-|x|^2/2}$
with respect to the Lebesgue measure on $\dR^n$ is the basic example of measure
$\mu $ satisfying~\eqref{eq-gross} with $C=2$. Since then, many results have presented measures
 satisfying such an inequality, among them the famous Bakry-Emery $\Gamma_2$ criterion, we refer
to Bakry \cite{bakry} and Ledoux \cite{ledoux} for further references and details on various applications of these inequalities.
\medskip

Let $\al\geq 1$ and define the probability measure $\mu_\al$ on $\R$ by
\begin{equation}
\label{eq-al}
\mu_\al(dx)=\frac{1}{Z_\al}e^{-\ABS{x}^\al}dx,
\end{equation}
where $Z_\al=\int e^{-\ABS{x}^\al}dx$. It is well-known that the probability
measure $\mu_\al$ satisfies a logarithmic Sobolev inequality~\eqref{eq-gross} if and only if
$\al\geq 2$. But for  $\al\in[1,2[$,  even if the measure $\mu_\al$ does not satisfy~\eqref{eq-gross},
it  satisfies a Poincar\'e inequality (or spectral gap inequality) which is  for every smooth enough function $f$,
\begin{equation}
\label{eq-poin}
\var{\mu_\al}{f}\leq C\int \ABS{\nabla f}^2d\mu_{\al},
\end{equation}
where $\var{\mu_\al}{f}=\int f^2 d\mu_\al-\PAR{\int f d\mu_\al}^2$ and $C<\infty$.

\medskip

Recall, see for example Section~1.2.6 of~\cite{logsob}, that if a probability measure on $\dR^n$ satisfies a
logarithmic Sobolev inequality with constant $C$ then it satisfies a Poincar\'e inequality with a constant less that
$C/2$.

\medskip

The problem  is then to interpolate between logarithmic Sobolev and  Poincar\'e inequalities, which will help us
to study further properties, such as concentration, of measures $\mu_\al^{\otimes n}$ for $\al\in[1,2]$ and $n\in\dN^*$.

A first answer was brought by Lata{\l}a-Oleszkiewicz in \cite{latala} and recently extended by
Barthe-Roberto in \cite{barthe-roberto}. Let $\mu$ be a probability measure on $\dR^n$, $\mu$ satisfies
 inequality $I_\mu(a)$ (for $a\in[0,1]$)  with constant $C>0$ if
for all $p\in[1,2[$,
\begin{equation}
\label{eq-lata}
\int f^2d\mu-\PAR{\int f^p d\mu}^{2/p}\leq C(2-p)^a\int \ABS{\nabla f}^2d\mu.
\end{equation}
A significant result of~\cite{latala} is that they prove that the measure  $\mu_\al^{\otimes^n}$
(for $\al\in[1,2]$, $n\in\dN^*$) satisfies such an inequality for a constant $C$ (independent of $n$)
 and with $a=2(\al-1)/\al$. And in~\cite{barthe-roberto} the authors
present a simple proof of the result of Lata{\l}a-Oleszkiewicz and describe the measures on the line
 which enjoy the same inequality.

\medskip

Our main purpose here will be to establish another type of interpolation between logarithmic Sobolev
and Poincar\'e inequalities, more directly linked to the structure of the usual logarithmic Sobolev
 inequalities, i.e. an inequality ``entropy-energy'' where we will modify the energy to enable us
to consider $\mu_\al$ measure. Note that this point of view was the one used by Bobkov-Ledoux
\cite{bobkov-ledoux2} when considering double sided exponential measure. Let us describe further
these modified logarithmic Sobolev inequalities.

Let $\al\in[1,2]$, $a>0$ and $\be\geq2$ satisfying  $1/\al+1/\be=1$, we note
$$
H_{a,\al}(x)= \left\{
\begin{array}{ll}
\disp\frac{x^2}{2} &\text{ if } \abs{x}\leq a \\
\disp a^{2-\be}\frac{\abs{x}}{\be}^\be+a^2\frac{\be-2}{2\be}\quad
&\text{ if } \abs{x}\geq a {\rm ~and~}\al\not=1\\
+\infty \quad&\text{ if } \abs{x}\geq a {\rm ~and~}\al=1.
\end{array}
\right.
$$
In Section~\ref{sec-def} we give definition and general properties of the following inequality
\begin{equation*}
\tag{$LSI_{{a,\al}}(C)$}
\ent{\mu}{f^2}\leq C \int H_{a,\al}\PAR{\frac{\nabla f}{f}}f^2 d\mu.
\end{equation*}
In particular we prove that inequality $LSI_{{a,\al}}$ satisfies some of the properties shared by Poincar\'e or Gross
logarithmic Sobolev inequalities (\eqref{eq-gross}
or~\eqref{eq-poin}), namely
tensorisation and perturbation. Note that in the case $\al=1$, it is exactly the inequality used by
Bobkov-Ledoux \cite{bobkov-ledoux2} and $\al=2$ is exactly the Gross logarithmic Sobolev inequality.

We present also a concentration property which is adapted to this inequality.
More precisely, if a measure $\mu$ satisfies the inequality $LSI_{{a,\al}}(C)$,
we have that if $f$ is a Lipschitz function on $\dR^n$ with $\NRM{f}_{Lip}\leq 1$ (with respect
to the Euclidean metric) then, there is $B>0$ such that for every $\la>0$ one has
\begin{equation}\label{superconc}
\mu_\al\PAR{f-\int fd\mu_\al\geq \la}\leq
\exp\PAR{-B\min\PAR{\la^\al,\la^2}}.
\end{equation}
This inequality was proved, for  $\al=1$,  by Maurey with the so called property $(\tau)$  and Bobkov-Ledoux in ~\cite{maurey,bobkov-ledoux2}.
Let us note that the cases $\al\geq 2$ are studied by Bobkov-Ledoux in~\cite{bobkov-ledoux1},
relying mainly on Brunn-Minkowski inequalities, and by Bobkov-Zegarlinski
in~\cite{bobkov-zeg} which refine the results presenting, via Hardy's inequality, some necessary and sufficient
condition for measures on the real line.
Let us note to finish that they use, for the case $\al\geq2$, $H_\be(x)=\ABS{x}^\be$ with $1/\al+1/\be=1$.

\medskip

In Section~\ref{subsec-tran},
we extend Otto-Villani's
theorem (see~\cite{villani})
for the relation with logarithmic Sobolev inequality and transportation inequality.
Let us define $L_{a,\al}$ by  $L_{a,\al}=H^*_{a,\al}$, the Legendre transform of $H_{a,\al}$.
We prove that if a probability measure $\mu$
on $\dR^n$ satisfies the inequality $LSI_{{a,\al}}(C)$ then there are $a'>0$ and $D>0$ such that
it satisfies also a transportation inequality: for all function $F$ on $\dR^n$, density of probability
with respect to $\mu$,
\begin{equation*}
\tag{$T_{{a',\al}}(D)$}
T_{L_{a',\al}}\PAR{F d\mu,d\mu}\leq D\ent{\mu}{F},
\end{equation*}
where
$$
T_{L_{a',\al}}\PAR{Fd\mu,d\mu}=\inf\BRA{\int L_{a',\al}\PAR{x-y}d\pi(x,y)},
$$
where the infimum is taken over the set of probabilities measures $\pi$ on $\dR^n\times\dR^n$ such that
$\pi$ has two margins $Fd\mu$ and $d\mu$. This inequality was introduced by Talagrand in~\cite{talagrand} for the case $\al=2$ and
$\al=1$. Let us note that the case $\al=1$ was also studied in~\cite{bgl} with exactly this form and  the case $\al\geq2$ was
studied
in~\cite{gentil-these}.

\medskip

In Section~\ref{sec-ex} we prove, as in~\cite{latala}, that the measure $\mu_\al$ defined in~\eqref{eq-al}
satisfies the inequality $LSI_{{a,\al}}(C)$. More precisely we
prove that there is $A,B>0$ such that $\mu_\al$ satisfies for all smooth function such that $f\geq0$ and $\int f^2d\mu_\al=1$,
\begin{equation*}
\ent{\mu_\al}{f^2}\leq A\var{\mu_\al}{f}+
B\int_{{f}\geq 2}\ABS{\frac{f'}{f}}^\be f^2d\mu_\al.
\end{equation*}
Due to the fact that $\mu_\al$ enjoys Poincar\'e inequality, $\mu_\al$  satisfies also
inequality $LSI_{{a,\al}}(C)$ for some constants $C>0$ and $a>0$.

Our method relies crucially on Hardy's inequality (see for example
\cite{logsob,bobkov-gotze,barthe-roberto}) we recall now: let $\mu,\nu$ be Borel measures on $\R^+$.
 Then the best constant $A$ so that every smooth function
$f$ satisfies
\begin{equation}
\label{eq-hardy}
\int_0^\infty \PAR{f(x)-f(0)}^2d\mu(x)\leq A\int_0^\infty f'^2d\nu
\end{equation}
is finite if and only if
\begin{equation}
\label{eq-condition}
B=\sup_{x>0} \mu\PAR{[x,\infty[} \int_0^x\PAR{\frac{d\nu^{ac}}{d\mu}}^{-1}dt
\end{equation}
is finite, where $\nu^{ac}$ is the absolutely continuous part of $\nu$ with respect to $\mu$.
Moreover, when $A$ is finite we have $$B\leq A\leq 4B.$$

\medskip
Finally in Section~\ref{section-autre} we will present some inequalities  satisfied by other measures.
More precisely, let $\phi$ be twice continuously differentiable and note the probability measure $\mu_\phi$ by,
\begin{equation}
\label{eq-mu}
\mu_{\phi}(dx)=\frac{1}{Z}e^{-\phi(x)}dx.
\end{equation}
Among them is considered
\begin{align*}
\phi(x) & =|x|^\al (\log(2+|x|))^\beta,\text{ with } \alpha\in]1,2[,\,\beta\in\R,
\end{align*}
which exhibits a modified logarithmic Sobolev inequality of function $H$ (different in nature from $H_{a,\al}$),
 and which is not covered by Lata{\l}a-Oleskiewickz inequality. We also present examples which are unbounded
perturbation of $\mu_\al$. We then derive new concentration inequalities in the spirit of Maurey \cite{maurey}
 or Bobkov-Ledoux \cite{bobkov-ledoux2}.

\section{Modified logarithmic Sobolev inequalities: definition and general properties}
\label{sec-def}

\subsection{Definitions and  classical properties}
\label{subsec-def}
Let $\al\in[1,2]$ and $\be\geq2$ satisfying  $1/\al+1/\be=1$ and let
$a>0$. Let define the functions $L_{a,\al}$ and  $H_{a,\al}$.

If $\al\in ]1,2]$ we note
$$
L_{a,\al}(x)= \left\{
\begin{array}{ll}
\disp\frac{x^2}{2} &\text{ if } \abs{x}\leq a \\
\disp a^{2-\al}\frac{\abs{x}}{\al}^\al+a^2\frac{\al-2}{2\al}\quad
&\text{ if } \abs{x}\geq a
\end{array}
\right.
$$
and
$$
H_{a,\al}(x)= \left\{
\begin{array}{ll}
\disp\frac{x^2}{2} &\text{ if } \abs{x}\leq a \\
\disp a^{2-\be}\frac{\abs{x}}{\be}^\be+a^2\frac{\be-2}{2\be}\quad
&\text{ if } \abs{x}\geq a
\end{array}
\right.
$$

If $\al=1$ we note
$$
L_{a,1}(x)= \left\{
\begin{array}{ll}
\disp\frac{x^2}{2} &\text{ if } \abs{x}\leq a \\
\disp a{\abs{x}}-\frac{a^2}{2}\quad
&\text{ if } \abs{x}\geq a
\end{array}
\right.
\text{ and }
H_{a,1}(x)= \left\{
\begin{array}{ll}
\disp\frac{x^2}{2} &\text{ if } \abs{x}\leq a \\
\disp \infty \quad&\text{ if } \abs{x}> a
\end{array}
\right.
$$

Let $n\in\dN^*$ and $x=(x_1,\cdots,x_n)\in\dR^n$, we note
$$
L_{a,\al}^{(n)}(x)=\sum_{i=1}^n L_{a,\al}(x_i) \text{ and }
H_{a,\al}^{(n)}(x)=\sum_{i=1}^n H_{a,\al}(x_i).$$
 Note that when there is no ambiguity we will drop the dependence in $n$ and note $L_{a,\al}$ instead of $L_{a,\al}^{(n)}$.

\medskip

Let us define the logarithmic Sobolev inequality of function  $H_{a,\al}$.

\begin{edefi}
Let $\mu$ be a probability measure on $\dR^n$,  $\mu$ satisfies a logarithmic Sobolev
inequality of function $H_{a,\al}$ with constant $C$, noted $LSI_{{a,\al}}(C)$, if for every $C^\infty$ and $L^2$ function $f$ on $\dR^n$
one has
\begin{equation}
\tag{$LSI_{{a,\al}}(C)$}
\ent{\mu}{f^2}\leq C \int H_{a,\al}\PAR{\frac{\nabla f}{f}}f^2 d\mu,
\end{equation}
where
$$\ent{\mu}{f^2}=\int f^2\log \frac{f^2}{\int f^2d\mu}d\mu\text{ and }
H_{a,\al}\PAR{\frac{\nabla f}{f}}=\sum_{i=1}^n H_{a,\al}\PAR{\frac{\partial f}{\partial x_i}\frac{1}{f}}.$$
It is supposed that $0/0=\infty$
\end{edefi}

We detail some properties of  $L_{a,\al}$ and
$H_{a,\al}$ in the following lemma.

\begin{elem}
\label{lemml}
 Functions $L_{a,\al}$ and
$H_{a,\al}$ satisfies:
\begin{enumerate}[i:]
\item If $\al\in]1,2]$, $L_{a,\al}$ and  $H_{a,\al}$ are $C^1$ on $\dR$.
\item $L_{a,\al}^*=H_{a,\al}$, where $L_{a,\al}^*$ is the Fenchel-Legendre transform
of   $L_{a,\al}$. Of course we have too $H_{a,\al}^*=L_{a,\al}$.

\item For all $t>0$ one has for all $x\in\dR$
$$
L_{a,\al}(tx)=t^2L_{\frac{a}{t},\al}(x),\quad H_{a,\al}(tx)=t^2H_{\frac{a}{t},\al}(x).
$$
\item Let $0\leq a\leq a'$, one has for all $x\in\dR^+$
$$
L_{a,\al}(x)\leq L_{a',\al}(x) ,\quad H_{a',\al}(x)\leq H_{a,\al}(x).
$$
\item If $\al\in]1,2]$, $L_{a,\al}$ and $H_{a,\al}$ are  strictly convex and
$$\lim_{\abs{x}\rightarrow\infty}\frac{H_{a,\al}(x)}{x}=\lim_{\abs{x}\rightarrow\infty}\frac{L_{a,\al}(x)}{x}=\infty.$$
\end{enumerate}
\end{elem}

The assumptions given on $\al$ and $\be$ are significant only for condition $iv$, and
condition $v$ is significant   for  Brenier-McCann-Gangbo's theorem, which is crucial
for the study of the link between modified logarithmic Sobolev inequalities and transportation inequalities of the next section.

Here are  some properties of the inequality $LSI_{{a,\al}}(C)$.

\begin{eprop}
\label{prop-facile}
\begin{enumerate}
\item This property is known under the name of tensorisation.

Let $\mu_1$ and $\mu_2$ two probability measures on
$\dR^{n_1}$ and $\dR^{n_2}$. Suppose that $\mu_1$ (resp. $\mu_2$) satisfies  the inequality $LSI_{{a,\al}}(C_1)$ (resp.
 $LSI_{{a,\al}}(C_2)$) then the probability $\mu_1\otimes\mu_2$ on $\dR^{n_1+n_2}$,
satisfies inequality $LSI_{{a,\al}}(D)$, where $D=\max\BRA{C_1,C_2}$.

\item  This property is known under the name of perturbation.

Let $\mu$ a measure on $\dR^n$ satisfying $LSI_{{a,\al}}(C)$. Let $h$ a bounded function on $\dR^n$ and defined
$\tilde{\mu}$ as
$$
d\tilde{\mu}=\frac{e^h}{{Z}}d\mu,
$$
where $Z=\int e^h d\mu$.

Then the measure $\tilde{\mu}$ satisfies the inequality  $LSI_{{a,\al}}(D)$ with $D=Ce^{\text{osc}(h)}$, where
$\text{osc}(h)=\sup(h)-\inf(h)$.
\item Link between  $LSI_{{a,\al}}(C)$ inequality with  Poincar\'e inequality.

Let $\mu$ a measure on  $\dR^n$. If $\mu$  satisfies $LSI_{{a,\al}}(C)$, then $\mu$ satisfies a Poincar\'e inequality
with the constant $C/2$. Let us recall that $\mu$ satisfies a Poincar\'e  inequality with constant $C$ if
\begin{equation}
\label{eq-poincare}
\var{\mu}{f}\leq C\int \ABS{\nabla f}^2d\mu,
\end{equation}
for all smooth function $f$.
\end{enumerate}

\end{eprop}

\begin{eproof}
One can find the details of the proof of the properties of tensorisation and  perturbation and the implication of the
Poincar\'e inequality in the chapter 1 and 3  of
\cite{logsob} (Section~1.2.6., Theorem~3.2.1 and Theorem~3.4.3).
\end{eproof}

\begin{erem}
We may of course define logarithmic Sobolev inequality of function $H$, where $H(x)$ is quadratic for
small values of $|x|$ and with convex, faster than quadratic, growth for large $|x|$. See Section~4
 for such examples. Note that Proposition \ref{prop-facile} is of course still valid for this kind of inequality.
These inequalities are also studied in a general case in \cite{ledoux} in Proposition~2.9.
\end{erem}

As in \cite{latala,bobkov-ledoux2}, by using the argument of Herbst,
one can give  precise estimates about concentration.

\begin{eprop}
  Assume that the probability measure $\mu$ on $\dR$ satisfies  the inequality $LSI_{{a,\al}}(C)$.
Let $F$ be Lipschitz function on $\dR$, then we get for $\la\geq0$,

$$
\mu(\ABS{F -\mu(F)}\geq \la)\leq
\left\{
\begin{array}{ll}
\disp2\exp\PAR{-K_\al\PAR{\la-aC\|F\|_{Lip}(2-\al)}^\al
-a^2\frac{2-\al}{2\al}}&\text{if }
\la\geq\frac{aC\|F\|_{Lip}}{2},\\
\disp2\exp\PAR{-\frac{2\la^2}{C\NRM{F}_{Lip}^2}}& \text{otherwise, }
\end{array}
\right.
$$
where $\displaystyle K_\al=\frac{ 2^\al(\al-1)^{1-\al}a^{2-\al}}
{\al C^{\al-1}\|F\|_{Lip}^\al}$.

Consider now $\mu^{\otimes n}$ and $F:\dR^n\rightarrow\dR$, $\mathcal{C}^1$, 
 such that
$\sum_{i=1}^n\ABS{\frac{\partial F}{\partial x_i}}^2\le 1$, then there
exists $\tilde K_\al$ (independent of $n$) such that
\begin{equation}\label{supconc}
\mu^{\otimes n}\PAR{\ABS{F -\mu^ {\otimes n}(F)}\geq \la}\leq
2\exp\left(-\tilde K_\al \min\left(\lambda^ 2,\lambda^\al\right)\right).
\end{equation}
\end{eprop}

\begin{eproof}
Assume that  $\int Fd\mu=0$.
Let us recall briefly Herbst's argument (see \cite{logsob} for more details). Denote $\Phi(t)=\int e^{t F}d\mu$,
and remark that $LSI_{a,\al}(C)$ applied to $f^2=e^{tF}$, using basic properties of $H_{a,\al}$, yields to
\begin{equation}\label{hehehe}
t\Phi'(t)-\Phi(t)\log \Phi(t)\le C H_{a,\al}\left(\frac{t\|F\|_{Lip}}{2}\right)\Phi(t)
\end{equation}
which, denoting $K(t)=(1/t)\log \Phi(t)$, entails
$$ K'(t)\le {C\over t^2}H_{a,\al}\left(\frac{t\|F\|_{Lip}}{2}\right).$$
Then, integrating, and using $K(0)=\int Fd\mu=0$, we obtain
\begin{equation}\label{laplace}
\Phi(t)\le \exp\left( Ct\int_0^t {1\over s^2}H_{a,\al}\left(\frac{s\|F\|_{Lip}}{2}\right)ds\right).\end{equation}
The Laplace transform of $F$ is then bounded by
$$\Phi(t)\leq
\left\{
\begin{array}{ll}
\disp\exp\PAR{Ct^\be\NRM{F}_{Lip}^\be\frac{a^{2-\be}}{2^\be\be(\be-1)}
+Ct a\NRM{F}_{Lip}\frac{\be-2}{2(\be-1)}
-Ca^2\frac{\be-2}{2\be}}& \text{ if }t\geq\frac{2a}{\NRM{F}_{Lip}},\\
\disp\exp\PAR{C\frac{\NRM{F}_{Lip}^2 t^2}{8}}&\text{ if }0\leq t\leq\frac{2a}{\NRM{F}_{Lip}}.
\end{array}
\right.
$$
For the $n$-dimensional extension, use the tensorisation property of $LSI_{a,\al}$ and 
$$
\sum_{i=1}^n H_{a,\al}\PAR{\frac{t}{2}\frac{\partial F}{\partial x_i}}\leq H_{a,\al}\PAR{\frac{t}{2}}.
$$ 
Then we can use the case of dimension 1.
\end{eproof}

\begin{erem}
For general logarithmic Sobolev of function $H$, we may obtain crude estimation of the concentration, at least for large
 $\lambda$. Indeed, using inequality (\ref{laplace}), we have directly that the concentration behavior is given by the
Fenchel-Legendre transform of $H$ for large values, see Section~\ref{section-autre}  for more details.
\end{erem}

\subsection{Link between inequality $LSI_{a,\al}(C)$ and transportation inequality}
\label{subsec-tran}

\begin{edefi}
Let $\mu$ be a probability measure on $\dR^n$,  $\mu$ satisfies a transportation inequality of function $L_{a,\al}$
with constant $C$, noted $T_{{a,\al}}(C)$, if for every function $F$, density of probability with respect to $\mu$,
one has
\begin{equation}
\tag{$T_{{a,\al}}(C)$}
T_{L_{a,\al}}\PAR{F d\mu,d\mu}\leq C\ent{\mu}{F},
\end{equation}
where
$$
T_{L_{a,\al}}\PAR{Fd\mu,\mu}=\inf\BRA{\int L_{a,\al}\PAR{x-y}d\pi(x,y)},
$$
where the infimum is taken over the set of probabilities measures $\pi$ on $\dR^n\times\dR^n$ such that
$\pi$ has two margins $Fd\mu$ and $\mu$.

\end{edefi}

Otto and Villani proved that a logarithmic Sobolev inequality implies a transportation inequality with a quadratic
cost (this is the case $\al=\be=2$), see \cite{villani,bgl}. With the notations of this paper they prove that
if $\mu$ satisfies the inequality $LSI_{\cdot,2}(C)$, (when $\al=2$
the constant $a$ is not any more a parameter in this case), then $\mu$ satisfies
the inequality  $T_{{\cdot,2}}(4C)$. In \cite{bgl} another case is studied, when
$\al=1$ and $\be=\infty$. In this first theorem we give an extension for the
other cases, where $\al\in[1,2]$.

\begin{ethm}
\label{thm-dir}
Let $\mu$ be a probability measure on $\dR^n$ and suppose that $\mu$ satisfies the inequality $LSI_{a,\al}(C)$.

Then $\mu$ satisfies the transportation inequality $T_{{\frac{aC}{2},\al}}(C/4)$.
\end{ethm}

\begin{eproof}
As in \cite{bgl}, we  use Hamilton-Jacobi equations. Let $f$ be a
Lipschitz bounded  function on $\dR^n$, and set
\begin{equation}
\label{eq-qt}
Q_t f(x)=\inf_{y\in\dR}\BRA{f(y)+tL_{\frac{aC}{2},\al}\PAR{\frac{x-y}{t}}},\,\,t>0,\,\,x\in\dR^n,
\end{equation}
and $Q_0 f=f$. The function $Q_t f$ is known as the Hopf-Lax solution of the Hamilton-Jacobi equation
\begin{equation*}
\left\{
\begin{array}{rl}
\disp
\frac{\partial v}{\partial t}(t,x)&=H_{\frac{aC}{2},\al}\PAR{\nabla v}(t,x),\,\,t>0,\,\,x\in\dR^n\\
\disp v(0,x)&=f(x),\,\,x\in\dR^n
\end{array}
\right.
\end{equation*}
see for example \cite{barles,evans}.

For $t\geq0$, define the function $\psi$ by
$$
\psi(t)=\int e^{\frac{4t}{C}Q_t f}d\mu.
$$
Since $f$ is Lipschitz and bounded  function one can prove that $Q_t f$ is also a  Lipschitz and bounded function on $t$
for almost every $x\in\dR^n$, then
$\psi$ is a $C^1$ function on $\dR^+$.
One gets
\begin{align*}
 \psi'(t)&=\int \frac{4}{C}Q_t fe^{\frac{4t}{C}Q_t
f}d\mu-\int\frac{4t}{C}H_{\frac{aC}{2},\al}\PAR{\nabla Q_t f}e^{\frac{4t}{C}Q_t
f}d\mu\\
&=\frac{1}{t}\ent{\mu}{e^{\frac{4t}{C}Q_t
f}}+\frac{1}{t}\psi(t)\log\psi(t)-\int\frac{4t}{C}H_{\frac{aC}{2},\al}\PAR{\nabla
Q_t f}e^{\frac{4t}{C}Q_t f}d\mu
\end{align*}
Let use inequality  $LSI_{a,\al}(C)$ to the function $\exp\left({\frac{2t}{C}Q_t f}\right)$ to get
\begin{equation*}
\psi'(t)\leq\frac{1}{t}\psi(t)\log\psi(t)
+\frac{C}{t}\PAR{\int
H_{a,\al}\PAR{\frac{2t}{C}\nabla Q_t f}e^{\frac{4t}{C}Q_t
f}d\mu-\int\frac{4t^2}{C^2}H_{\frac{aC}{2},\al}\PAR{\nabla Q_t
f}e^{\frac{4t}{C}Q_t f}d\mu}.
\end{equation*}
Due to the property of $H_{a,\al}$ (see Lemma~\ref{lemml}),
$$
H_{a,\al}\PAR{\frac{2t}{C}\nabla Q_t
f}=\frac{4t^2}{C^2}H_{\frac{aC}{2t},\al}\PAR{\nabla Q_t f}.
$$
Then for all  $t\in[0,1]$, one has
$$
H_{a,\al}\PAR{\frac{2t}{C}\nabla Q_t f}\leq
\frac{4t^2}{C^2}H_{\frac{aC}{2},\al}\PAR{\nabla Q_t f}.
$$
Then
$$
\forall t\in[0,1],\quad t\psi'(t)-\psi(t)\log\psi(t)\leq 0
$$
After integration on $[0,1]$, we have
$$
\psi(1)\leq \exp\frac{\psi'(0)}{\psi(0)},
$$
from where
\begin{equation}
\label{eq-bg}
\int e^{\frac{4}{C}Q_1 f}d\mu\leq e^{\int\frac{4}{C}f d\mu}.
\end{equation}

Since
$$\ent{\mu}{F}=\sup\BRA{\int Fgd\mu,\,\,\int e^gd\mu\leq1},$$
we have with $g=\frac{4}{C}Q_1 f-\int\frac{4}{C}fd\mu$,
$$
\int F\PAR{Q_1 f-\int fd\mu}d\mu\leq \frac{C}{4}\ent{\mu}{F}.
$$

Let take the supremum on the set of Lipschitz function $f$, the Kantorovich-Rubinstein's theorem
applied to  the distance $T_{L_{a,\al}}(Fd\mu,d\mu)$, see \cite{villani03},
implies that
$$
T_{L_{\frac{aC}{2}},\al}(F d\mu,d\mu)\leq \frac{C}{4}\ent{\mu}{F}.
$$
\end{eproof}

As it is also the case in quadratic case, when the measure is  log-concave
one can prove that a transportation inequality implies a logarithmic
Sobolev inequality.

\medskip

{ \it For the next theorem we suppose that the function of transport given
by the theorem of Brenier-Gangbo-McCann is a $\mathcal{C}^2$ function. Such a regularity result is
outside the scope of this paper and we refer to Villani \cite{villani03} for further discussions
 around this problem. However we show here, that once this result assumed, the methodology presented
 in Bobkov-Gentil-Ledoux \cite{bgl}, for the exponential measure,  still works. }

 \begin{ethm}
\label{thm-inv}
Let $\mu$ be a probability measure on $\dR^n$. Assume that
 $$\mu(dx)=e^{-\phi(x)}dx$$ where  $\phi$ is a convex function on $\dR^n$.

If  $\mu$ satisfies the inequality $T_{{a,\al}}(C)$ then for all $\la>C$,
$\mu$ satisfies the logarithmic Sobolev inequality $LSI_{{\frac{a}{2\la},\al}}\PAR{\frac{4\la^2}{\la-C}}$.
 \end{ethm}
\begin{eproof}
Let note $F$ density of probability with respect to
$\mu$. Assume that $F$ is $\mathcal{C}^2$, the general case can result by density.

By the  Brenier-Gangbo-McCann's theorem, see
\cite{brenier91,gangbo-mccann}, there exists a function
$\Phi$ such that
$$
S=\text{Id}-{\nabla H_{a,\al}}\circ\nabla\Phi,
$$
transports $Fd\mu$ to the measure $\mu$, for every measurable bounded function $g$
$$
\int g(S) Fd\mu=\int g d\mu.
$$
The function $\Phi$ is a $L_{a,\al}$-concave function and  if $\Phi$ is $\mathcal{C}^2$,
a classical argument of convexity (see chapter 2 of \cite{villani03}),
one has
$D\left[{\nabla H_{a,\al}\circ\nabla\Phi(x)}\right]$ is diagonalizable with real eigenvalues, all less than~1.

According to the assumption made on function $\Phi$,  one can assume that $S$
is sufficiently smooth and  we obtain for $x\in\dR^n$,
\begin{equation}
\label{eq-monge}
F(x)e^{-\phi(x)}=e^{-\phi\circ S(x)}\det\PAR{\nabla S(x)}.
\end{equation}

Moreover this function gives the optimal transport, i.e.
$$
T_{L_{a,\al}}\PAR{F    d\mu,d\mu}=
\int L_{a,\al}\PAR{{\nabla {H_{a,\al}}}\circ\nabla\Phi}Fd\mu.
$$

Then by~\eqref{eq-monge}, one has for $x\in\dR^n$,
\begin{align*}
\log F(x) =\phi(x)-\phi(x-{\nabla H_{a,\al}}\circ\nabla\Phi(x))+
\log{\det\PAR{\text{Id}-D\left[{\nabla H_{a,\al}\circ\nabla\Phi(x)}\right]}}.
\end{align*}

Then since $D\left[{\nabla H_{a,\al}\circ\nabla\Phi(x)}\right]$ is diagonalizable with real eigenvalues, all less than~1,
we get
$$
\log{\det\PAR{\text{Id}-D\left[{\nabla H_{a,\al}\circ\nabla\Phi(x)}\right]}}\leq
-\text{div}\PAR{{\nabla H_{a,\al}}\circ\nabla\Phi(x)}.
$$

Since $\phi$ is convex we have $\phi(x)-\phi(x-{\nabla H_{a,\al}}\circ\nabla\Phi(x))\leq
{\nabla H_{a,\al}}\circ\nabla\Phi(x)\cdot\nabla\phi(x)$ and we obtain
\begin{align*}
\ent{\mu}{F}\leq \int\BRA{{\nabla H_{a,\al}}\circ\nabla\Phi(x)\cdot\nabla\phi(x)
-\text{div}\PAR{{\nabla H_{a,\al}}\circ\nabla\Phi(x)}}F(x)d\mu(x),
\end{align*}
after integration  by parts
$$
\ent{\mu}{F}\leq\int\nabla F\cdot {\nabla H_{a,\al}}\circ\nabla\Phi d\mu.
$$

Let $\la>0$ and let use Young inequality for the combined functions $L_{a,\al}$ and $H_{a,\al}$
$$
\la \frac{\nabla F}{F} \cdot{\nabla H_{a,\al}}\circ\nabla\Phi\leq
H_{a,\al}\PAR{\la\frac{\nabla F}{ F}}+L_{a,\al}\PAR{{\nabla H_{a,\al}}\circ\nabla\Phi}.
$$
Thus
\begin{align*}
\ent{\mu}{F}&\leq\frac{1}{\la}\int H_{a,\al}\PAR{\la\frac{\nabla F}{F}} F d\mu
+\frac{1}{\la}\int L_{a,\al}\PAR{{\nabla H_{a,\al}}\circ\nabla\Phi} Fd\mu\\
&\leq {\la} \int H_{\frac{a}{\la},\al}\PAR{\frac{\nabla F}{F}} F d\mu
+\frac{1}{\la} T_{L_{a,\al}}\PAR{F d\mu,d\mu}.
\end{align*}

Thus if  $\mu$ satisfies the inequality $T_{{a,\al}}(C)$ we get for all $\la>C$
$$
\ent{\mu}{F}\leq\frac{\la^2}{\la-C} \int H_{\frac{a}{\la},\al}\PAR{\frac{\nabla F}{F}}F d\mu.
$$
Let us note now $f^2=F$, we get
\begin{align*}
\ent{\mu}{f^2}&\leq\frac{\la^2}{\la-C} \int H_{\frac{a}{\la},\al}\PAR{2\frac{\nabla f}{f}}f^2 d\mu\\
&\leq\frac{4\la^2}{\la-C} \int H_{\frac{a}{2\la},\al}\PAR{\frac{\nabla f}{f}}f^2 d\mu.
\end{align*}
Then $\mu$ satisfies, for all $\la>C$  inequality  $LSI_{\frac{a}{2\la},\al}\PAR{\frac{4\la^2}{\la-C}}$.
\end{eproof}

\begin{erem}
\label{essai}
One can summarizes Theorem~\ref{thm-dir} and \ref{thm-inv} by the 
following diagram (under assumption of Theorem~\ref{thm-inv}):
\begin{align*}
{LSI_{a,\al}(C)}&\rightarrow T_{\frac{aC}{2},\al}(C/4)\\
{T_{a,\al}(C)}&\rightarrow \BRA{LSI_{\frac{a}{2\la},\al}\PAR{\frac{4\la^2}{\la-C}}}_{\la>C}.
\end{align*}
Notice, as it is the case for the traditional logarithmic Sobolev inequality,
 than there is a loss on the level of the constants in the direction transportation inequality implies
 logarithmic Sobolev inequality. When $\al=\be=2$, we get as in \cite{villani},
${T_{\cdot,2}(C)}\rightarrow {LSI_{\cdot,2}(16C)}$. As in~\cite{villani}, Theorem~\ref{thm-inv} can be modified in the case
$\text{Hess}( \phi)\geq \la \text{Id}$, where $\la\in\dR$.
\medskip

Also let us notice that as in the quadratic  case we do not know if these two inequalities are equivalent.
\end{erem}

As in Proposition~\ref{prop-facile}, here are  some properties of the inequality $T_{L_{a,\al}}(C)$.

\begin{eprop}

\begin{enumerate}
\item Let us recall Marton's theorem on concentration inequality.

Assume that $\mu$ satisfies a transportation inequality  $T_{L_{a,\al}}(C)$ then $\mu$ satisfies the following
 concentration inequality
\begin{equation*}
\label{eq-marton}
\forall A\subset R^n,\quad \text{with} \quad \mu(A)\geq
\frac{1}{2},\quad \mu(\PAR{A_r}^c)\leq
2e^\PAR{-\frac{1}{C}L_{a,\al}(r)},
\end{equation*}
where $\PAR{A_r}^c=\BRA{x\in\dR^n,\,\,d(A,x)\geq r}$.

\item As  in  Proposition~\ref{prop-facile},  the
properties of tensorisation are also valid for transportation inequality $T_{a,\al}(C)$.

Let $\mu_1$ and $\mu_2$ be two probability measures on
$\dR^{n_1}$ and $\dR^{n_2}$. Suppose that $\mu_1$ (resp. $\mu_2$) satisfies
the inequality $T_{{a,\al}}(C_1)$ (resp.
 $T_{{a,\al}}(C_2)$) then the probability $\mu_1\otimes\mu_2$ on $\dR^{n_1+n_2}$,
satisfies inequality $T_{{a,\al}}(D)$, where $D=\max\BRA{C_1,C_2}$.

\item If the measure $\mu$ verifies $T_{{a,\al}}(C)$,
then $\mu$ satisfies a Poincar\'e inequality~\eqref{eq-poincare}
with the constant~$C$.
\end{enumerate}
\end{eprop}

\begin{eproof}
 The demonstration of $i$, $ii$
of these results is a simple adaptation of the traditional case,
we return to the references for proofs (for example chapters 3, 7 and 8 of~\cite{logsob}).

The proof of $iii$ is an adaptation of the quadratic case. Suppose that $\mu$ satisfies a $T_{{a,\al}}(C)$.
By a classical argument of Bobkov-G\"otze, the measure $\mu$ satisfies the dual form of $T_{{a,\al}}(C)$ which is
the inequality~\eqref{eq-bg},
\begin{equation}
\label{eq-bg2}
\int e^{\frac{1}{C}Q_1 f}d\mu\leq e^{\int\frac{1}{C}f d\mu},
\end{equation}
where $Q_1 f$ is defined as in~\eqref{eq-qt} with the function $L_{a,\al}$.

Let note $f=\ep g$ with $g$, $\mathcal{C}^1$ and bounded, we get
\begin{align*}
Q_1 f(x)=Q_1 (\ep g)(x)&=\ep\inf_{z\in\dR^n}\BRA{g(x-\ep z)+\ep L_{\frac{a}{\ep},\al}(z)}\\
& =\ep {g(x)-\frac{\ep^2}{2}\ABS{\nabla g}^2+o(\ep^2)}
\end{align*}
Then we obtain by~\eqref{eq-bg2},
$$
1+\frac{\ep}{C}\int gd\mu-\frac{\ep^2}{2C}\int \ABS{\nabla g}^2d\mu+
\frac{\ep^2}{2C^2}{\int g^2d\mu}\leq1+\frac{\ep}{C}\int gd\mu
+\frac{\ep^2}{2 C^2}\PAR{\int gd\mu}^2+o(\ep^2),
$$
 imply that
$$
\var{\mu}{g}\leq C\int \ABS{\nabla g}^2d\mu.
$$
\end{eproof}

Unfortunately, as in the traditional case of the transportation inequality,
we do not know if this one has property of perturbation as for
inequality  $LSI_{a,\al}(C)$.
\section{An important example on $\dR$,  the measure   $\mu_\al$}
\label{sec-ex}

Let $\al\geq1$ and  define the probability measure $\mu_\al$ on $\R$ by
$$\mu_\al(dx)=\frac{1}{Z_\al}e^{-\ABS{x}^\al}dx,$$
where $Z_\al=\int e^{-\ABS{x}^\al}dx$.

\begin{ethm}
\label{thm-the}
Let $\al\in]1,2]$. There exists $A,B>0$ such  that the measure $\mu_\al$ satisfies the following
modified logarithmic Sobolev inequality, for any smooth function $f$ on $\R$ such that $f\geq0$ and $\int f^2d\mu_\al=1$  we have
\begin{equation}
\label{eq-hardy2}
\ent{\mu_\al}{f^2}\leq A\var{\mu_\al}{f}+
B\int_{{f}\geq 2}\ABS{\frac{f'}{f}}^\be f^2d\mu_\al,
\end{equation}
where $1/\al+1/\be=1$ and
$$
\ent{\mu_\al}{f^2}=\int f^2\log \frac{f^2}{\int f^2d\mu_\al}d\mu_\al\text{ and }
\var{\mu_\al}{f}=\int f^2d\mu_\al-\PAR{\int fd\mu_\al}^2.
$$
In the extreme case, $\al=1$,  we obtain the following inequality: for all $f$ smooth enough such that $|f'|\le 1$,
\begin{equation}
\label{eq-extrem}
\ent{\mu_\al}{f^2}\leq A\var{\mu_\al}{f}.
\end{equation}
\end{ethm}

\begin{ecor}
\label{cor-the}
Assume that $f$ is a smooth function on $\dR$. Then we obtain the following estimation
\begin{equation}
\label{eq-cor}
\ent{\mu_\al}{f^2}\leq A\var{\mu_\al}{f}+
B\int_{\Omega}\ABS{\frac{f'}{f}}^\be f^2d\mu_\al,
\end{equation}
where
$$
\Omega=\BRA{f_+\geq2\sqrt{\int f_+^2d\mu_\al}}\cup\BRA{f_-\geq2\sqrt{\int f_-^2d\mu_\al}},
$$
$f_+=\max(f,0)$ and $f_-=\max(-f,0)$.
\end{ecor}

\begin{eproof}
We have
$f^2=f_+^2+f_-^2$. Then
\begin{eqnarray*}
\ent{\mu_\al}{f^2}& = &\sup\BRA{\int f^2gd\mu_\al \text{ with }\int e^g d\mu_\al\leq1}\\
& = &  \sup\BRA{\int f_+^2gd\mu_\al+\int f_-^2gd\mu_\al \text{ with }\int e^g d\mu_\al\leq1}\\
&\leq & \ent{\mu_\al}{f_+^2}+\ent{\mu_\al}{f_-^2}.
\end{eqnarray*}

By  Theorem~\ref{thm-the} there
 exists $A,B>0$ independent of $f$ such that
\begin{equation*}
\label{eq-++}
\ent{\mu_\al}{f_+^2}\leq A\var{\mu_\al}{f_+}+
B\int_{\Omega_+}\ABS{\frac{f_+'}{f_+}}^\be f_+^2d\mu_\al,
\end{equation*}
\begin{equation*}
\ent{\mu_\al}{f_-^2}\leq A\var{\mu_\al}{f_-}+
B\int_{\Omega_-}\ABS{\frac{f_-'}{f_-}}^\be f_-^2d\mu_\al,
\end{equation*}
where $\Omega_+=\BRA{f_+\geq2\sqrt{\int f_+^2d\mu_\al}}$ and  $\Omega_-=\BRA{f_-\geq2\sqrt{\int f_-^2d\mu_\al}}$.

To conclude, it is enough to notice that
\begin{eqnarray*}
\var{\mu_\al}{f_+}+\var{\mu_\al}{f_-}& = &\int f^2d\mu_\al-\PAR{\PAR{\int f_+d\mu_\al}^2
+\PAR{\int f_-d\mu_\al}^2}\\
&\leq & \var{\mu_\al}{f},
\end{eqnarray*}
and
\begin{eqnarray*}
\int_{\Omega_+}\ABS{\frac{f_+'}{f_+}}^\be f_+^2d\mu_\al+
\int_{\Omega_-}\ABS{\frac{f_-'}{f_-}}^\be f_-^2d\mu_\al=
\int_{\Omega}\ABS{\frac{f'}{f}}^\be f^2d\mu_\al.
\end{eqnarray*}
\end{eproof}

It implies there are $a_\al>0$ and $C_\al<\infty$, such that  $\mu_\al$ satisfies a logarithmic Sobolev 
inequality of function $H_{a_\al,\al}$ with constant $C_\al$.

Indeed, this is clear that $\mu_\al$ satisfies a
Poincar\'e inequality, (see chapter~6 of \cite{logsob}),
with constant~$\la_\al<\infty$,
$$
\var{\mu_\al}{f}\leq\la_\al\int f'^2d\mu_\al.
$$
Then, by inequality~\eqref{eq-hardy2}, we obtain for any smooth function $f$ on $\dR$,
$$
\ent{\mu_\al}{f^2}\leq A\la_\al\int{{{f'}^2}d\mu_\al}+
B\int\ABS{\frac{f'}{f}}^\be f^2d\mu_\al.
$$

\bigskip

Let us give a few hint on the proof of the Theorem~\ref{thm-the}, which will enable us to present key auxiliary lemmas.
 We first use the
 following inequality

\begin{equation}
\int f^2\ln f^2d\mu_\al  \leq 5\int (f-1)^2d\mu_\al +\int (f-2)_+^2\ln (f-2)_+^2d\mu_\al
\label{triv-ineq}
\end{equation}
where it is obvious that truncation arguments are crucial. We will then need the following lemma:

\begin{elem}
\label{lem-var}
Let $\mu$ be a probability measure on $\dR$ and let $f\geq0$ such that $\int f^2d\mu=1$ then we obtain
\begin{enumerate}[i:]
\item $\displaystyle \int (f-1)^2\leq 2\var{\mu}{f}.$
\item $\displaystyle \int_{f\ge2}f^2d\mu \leq 8 \var{\mu}{f}.$
\item $\displaystyle\int_{f\geq 2}f^2\ln f^2d\mu\leq \frac{\ln 4}{\ln 4 -1} \ent{\mu}{f^2}<4\ent{\mu}{f^2}.$
\end{enumerate}
\end{elem}

\begin{eproof}
{\it i}. We have
$$
\int (f-1)^2d\mu=\var{\mu}{f}+\PAR{1-\int f d\mu}^2.
$$
$\int f^2d\mu=1$ imply that $0\leq \int f d\mu\leq 1$, then $\PAR{\int f d\mu}^2\leq \int f d\mu$. Then
$$
\int (f-1)^2d\mu\leq\var{\mu}{f}+(\var{\mu}{f})^2,
$$
but since $\int f^2d\mu=1$, $\var{\mu}{f}\leq 1$,  then
$$
\int (f-1)^2d\mu\leq 2\var{\mu}{f}.
$$

{\it ii}. One verifies trivially that when $x\ge2$, $x^2\le {4}(x-1)^2$ and apply {\it i}.

{\it iii}. Let us give the proof given in \cite{gui-ca}.

If $x>0$ we have $x\ln x+1-x\geq 0$ which yields
$$
\int_{f\leq 2}f^2\ln f^2d\mu+\mu\PAR{f\leq 2}-\int_{f\leq 2}f^2d\mu\geq 0,
$$
hence
\begin{eqnarray*}
\ent{\mu}{f^2}\geq\int_{f\geq 2} f^2\ln f^2d\mu-\int_{f\geq 2}f^2d\mu.
\end{eqnarray*}
Since
$$\int_{f\geq 2}f^2d\mu\leq\frac{1}{\ln 4}\int_{f\geq 2} f^2\ln f^2d\mu,$$
we obtain
\begin{eqnarray*}
\ent{\mu}{f^2}\geq\PAR{1-\frac{1}{\ln 4}}\int_{f\geq 2} f^2\ln f^2d\mu.
\end{eqnarray*}
\end{eproof}

Recall the Hardy's inequality presented in the introduction: let $\mu,\nu$ be Borel measures on $\R^+$,
the best constant $A$ so that every smooth function
$f$ satisfies
\begin{equation}
\label{eq-hardybis}
\int_0^\infty \PAR{f(x)-f(0)}^2d\mu(x)\leq A\int_0^\infty f'^2d\nu
\end{equation}
is finite if and only if
\begin{equation}
\label{eq-conditionbis}
B=\sup_{x>0} \mu\PAR{[x,\infty[} \int_0^x\PAR{\frac{d\nu^{ac}}{d\mu}}^{-1}dt
\end{equation}
is finite and when $A$ is finite we have $$B\leq A\leq 4B.$$

We then present different proof of the desired inequality, starting from (\ref{triv-ineq}),
according to the value of $\ent{\mu}{f^2}$, in which Hardy's inequality plays a crucial role.
First, when the entropy is large we will need

\medskip

\begin{elem}
\label{lem-h}
Let $h$ defined as follow.

\begin{equation*}
h(x)=\left\{
\begin{array}{lr}
1&\text{if}\,\,\ABS{x}\leq1\\
\ABS{x}^{2-\al}&\text{if}\,\,\ABS{x}\geq1
\end{array}
\right.
\end{equation*}

Then
there exists ${C_h}>0$ such that  for every smooth function $g$ we have
\begin{equation}
\label{eq-lemh}
\ent{\mu_\al}{g^2}\leq {C_h}\int g'^2 hd\mu_\al.
\end{equation}
\end{elem}

\begin{eproof}
We use Theorem~3 of \cite{barthe-roberto} which is a refinement of the criterion of a Bobkov-G\"otze theorem
(see Theorem~5.3 of \cite{bobkov-gotze}).

The constant ${C_h}$ satisfies $\max(b_{-},b_{+})\leq C\leq \max(B_{-},B_{+})$ where
$$
b_{+}=\sup_{x>0} \mu_\al([x,+\infty[)\log\PAR{1+\frac{1}{2\mu_\al([x,+\infty[)}}
\int_0^x  {Z_\al}\frac{e^{\ABS{t}^\al}}{h(t)}dt,
$$

$$
b_{-}=\sup_{x<0} \mu_\al(]-\infty,x])\log\PAR{1+\frac{1}{2\mu_\al(]-\infty,x])}}
\int_x^0  {Z_\al}\frac{e^{\ABS{t}^\al}}{h(t)}dt,
$$

$$
B_{+}=\sup_{x>0} \mu_\al([x,+\infty[)\log\PAR{1+\frac{e^2}{\mu_\al([x,+\infty[)}}
\int_0^x  {Z_\al}\frac{e^{\ABS{t}^\al}}{h(t)}dt,
$$

$$
B_{-}=\sup_{x<0} \mu_\al(]-\infty,x])\log\PAR{1+\frac{e^2}{\mu_\al([-\infty,x[)}}
\int_x^0  {Z_\al}\frac{e^{\ABS{t}^\al}}{h(t)}dt.
$$

An easy approximation prove that for large positive $x$
\begin{equation}
\label{eq-approx}
\mu_\al([x,\infty[)=\int_x^\infty \frac{1}{Z_\al}e^{-\ABS{t}^\al}dt
\sim_\infty \frac{1}{Z_\al\al x^{\al-1}}e^{-x^\al}
\end{equation}
$$
\int_0^x  {Z_\al}\frac{e^{\ABS{t}^\al}}{h(t)}dt\sim_\infty \frac{Z_\al}{x}e^{x^\al},
$$
and one may prove same equivalent for negative $x$. A simple calculation then yields that
constants $b_{+}$, $b_{-}$, $B_{+}$ and  $B_{-}$ are finite and the lemma is proved.
\end{eproof}

Note that the function $h$ is the smallest function
such that the constant ${C_h}$ in the inequality~\eqref{eq-lemh} is finite,
it ``saturates'' the inequality on infinity.

\medskip

In the case of small entropy, we will use so-called $\Phi$-Sobolev inequalities (even if our context
is less general), see Chafa{\"\i} \cite{chafai} for a comprehensive review, and Barthe-Cattiaux-Roberto
 \cite{ca-ba-ro} for a general approach in the case of measure $\mu_\al$.

\begin{elem}
\label{lem-g}
Let $g$ be defined on $[T,\infty[$ with $T\in[T_1,T_2[$ for some fixed $T_1,T_2$,
$$g(T)=2,\,\, g\geq2\,\, \text{and}\,\, \int_{T}^\infty g^2d\mu_\al\leq13,$$ then
\begin{equation}
\label{eq-lemme}
\int_{T}^\infty (g-2)_+^2\Phi(g^2)\mu_\al\leq
 {C_g}\int_{[T,\infty[\cap\BRA{g\geq2}} {g'}^2d\mu_\al,
\end{equation}
where $\Phi(x)=\ln^{\frac{2(\al-1)}{\al}}(x)$. The constant ${C_g}$ depend on the measure $\mu_\al$ but does not depend
on the value of $T\in[T_1,T_2]$.
\end{elem}

\begin{eproof}
Let use Hardy's inequality as explained in the introduction.
We have $g(T)=2$. We apply inequality~\eqref{eq-hardy} with the function $(g-2)_+$
and the following measures
$$d\mu=\PAR{\ln g^2}^{\frac{2(\al-1)}{\al}}d\mu_\al\,\,\text{and}\,\,\nu=\mu_\al.$$
Then  the constant $C$ in inequality~\eqref{eq-lemme} is finite if and only if
$$
B=\sup_{x>T}\int_{T}^x Z_\al e^{\ABS{t}^\al}dt\int_x^\infty\PAR{\ln g^2}^{\frac{2(\al-1)}{\al}}d\mu_\al,
$$
is finite.

Since ${2(\al-1)}/{\al}<1$  the function
$x\rightarrow\PAR{\ln x}^{\frac{2(\al-1)}{\al}}$ is concave on $[4,\infty[$.
By Jensen inequality we obtain for all $x\geq T$,
$$
\int_x^\infty\PAR{\ln g^2}^{\frac{2(\al-1)}{\al}}d\mu_\al
\leq  \ln^{\frac{2(\al-1)}{\al}}\PAR{\frac{\int_{x}^\infty g^2d\mu_\al}
{\mu_{\al}\PAR{[x,\infty[}}}\mu_\al\PAR{[x,\infty[}.
$$
Then by the property of $g$ we have
\begin{eqnarray*}
B & \leq & \sup_{x>T}\int_{T}^x Z_\al e^{\ABS{t}^\al}dt
\ln^{\frac{2(\al-1)}{\al}}\PAR{\frac{13}{\mu_\al\PAR{[x,\infty[}}}\mu_\al\PAR{[x,\infty[}\\
& \leq & \sup_{x>T_1}\int_{T_1}^x Z_\al e^{\ABS{t}^\al}dt
\ln^{\frac{2(\al-1)}{\al}}\PAR{\frac{13}{\mu_\al\PAR{[x,\infty[}}}\mu_\al\PAR{[x,\infty[}.
\end{eqnarray*}
Using the approximation given in equality~\eqref{eq-approx} we prove that $B$ is finite, bounded
by a constant ${C_g}$ which does not depend on $T$.
\end{eproof}

As said before, we divide the proof of Theorem \ref{thm-the} in two parts: large and small entropy,
both in the case of positive function.  Let us now
 present the proof in the case of large entropy.

\medskip
{\it Large entropy case}.


\begin{eprop}
\label{prop-grand2}
Suppose that $\al\in]1,2]$.
There exists $A,B>0$ such that for any functions $f\geq0$ satisfying
$$\int f^2 d\mu_\al=1\,\,\text{and}\,\,\ent{\mu_\al}{f^2}\geq1$$
 we have
\begin{equation}
\label{eq-grand}
\ent{\mu_\al}{f^2}\leq
A\var{\mu_\al}{f}+B\int_{f\geq2}\ABS{\frac{f'}{f}}^\be f^2d\mu_\al.
\end{equation}
If $\al=1$, when $|f'|\le 1$, then $\ent{\mu_\al}{f^2}\leq A\var{\mu_\al}{f}$.
\end{eprop}

\bigskip

{\noindent {\emph{\textbf{Proof of Proposition~\ref{prop-grand2}}}}\\~}
\proofbegin~
Let $f\geq0$ satisfying $\int f^2d\mu_\al=1$.

A careful study of the function
$$
x\rightarrow -x^2\ln x^2 +5(x-1)^2+x^2-1+(x-2)_+^2\ln (x-2)_+^2
$$
proves that for every $x\in \R$
$$
x^2\ln x^2\leq 5(x-1)^2+x^2-1+(x-2)_+^2\ln (x-2)_+^2.
$$

Then we obtain by Lemma~\ref{lem-var}.i, recalling that $\int f^2d\mu_\al=1$ and $f\geq0$,
\begin{eqnarray*}
\int f^2\ln f^2d\mu_\al & \leq &5\int (f-1)^2d\mu_\al +\int (f^2-1)d\mu_\al+
\int (f-2)_+^2\ln (f-2)_+^2d\mu_\al\\
 & \leq &10\var{\mu_\al}{f}+ \int (f-2)_+^2\ln (f-2)_+^2d\mu_\al
\end{eqnarray*}

which is the announced starting point inequality (\ref{triv-ineq}).

Since $\int f^2 d\mu_\al=1$, one can easily prove that
$$
\int (f-2)_+^2d\mu_\al\leq 1,
$$
then $\int (f-2)_+^2\ln (f-2)_+^2d\mu_\al\leq \ent{\mu_\al}{(f-2)_+^2},$
and
\begin{eqnarray*}
\ent{\mu_\al}{f^2}\leq &10\var{\mu_\al}{f}+\ent{\mu_\al}{(f-2)_+^2}.
\end{eqnarray*}

Hardy's inequality of Lemma~\ref{lem-h} with $g=(f-2)_+$ gives

\begin{equation}
\label{eq-br}
\ent{\mu_\al}{(f-2)_+^2}\leq {C_h}\int (f-2)_+'^2 hd\mu_\al
={C_h}\int_{f\geq2} f'^2 hd\mu_\al.
\end{equation}

For $p\in]1,2[$ and $q>2$ such that  and $1/p+1/q=1$ and we have for every $x,y>0$ by Young inequality,
\begin{equation}
\label{eq-young}
xy\leq \frac{x^p}{p}+\frac{y^q}{q}.
\end{equation}

If $\al=1$, then if $|f'|\le1$, then there exists $C>0$ such that
$${C_h}\int_{f\geq2} f'^2 hd\mu_\al\le  C\var{\mu_\al}{f}+{1\over 2}\ent{\mu_\al}{f^2}$$
where we used Lemma \ref{lem-var}.ii and the large entropy case. We then deduce the result when $\al=1$.

Consider then $\al\in]1,2]$ and $\be=\al/(\al-1)$. Let $p=\be/2$ and $q=\be/(\be-2)$. Let $\e>0$ and let apply
inequality~\eqref{eq-young} to the right term of~\eqref{eq-br}, we obtain
$$
\frac{1}{\e^{(\be-2)/\be}}\PAR{\frac{f'}{f}}^2 \e^{(\be-2)/\be}h\leq
\frac{2}{\be\e^{(\be-2)/2}}{\ABS{\frac{f'}{f}}}^\be+\frac{\be-2}{\be}\e h^{\be/(\be-2)},
$$
then
\begin{equation*}
\ent{\mu_\al}{{(f-2)_+}^2}  \leq
\frac{2{C_h}}{\be\e^{(\be-2)/2}}
\int_{f\geq2}\ABS{\frac{f'}{f}}^\be f^2d\mu_\al+
\frac{\be-2}{\be}{C_h}\e\int_{f\geq2} h^{\be/(\be-2)}f^2d\mu_\al
\end{equation*}

Let $\mu$ a probability measure, then we have for every function $f$ such that $\int f^2d\mu=1$,
and for every measurable function $g$
$$
\int f^2g d\mu \leq \ent{\mu}{f^2} + \log{\int e^g d\mu}.
$$
Let $\eta>0$ and  we apply the previous inequality with $g=\eta h^{\be/(\be-2)}$,
\begin{multline}
\ent{\mu_\al}{{(f-2)_+}^2}\leq \\
\frac{2{C_h}}{\be\e^{(\be-2)/2}}
\int_{f\geq2}\ABS{\frac{f'}{f}}^\be f^2d\mu_\al +
\frac{(\be-2) {C_h}\e}{\be\eta}
\PAR{ \ent{\mu_\al}{f^2}+ \log{\int \exp\PAR{\eta h^{\be/(\be-2)}}d\mu_\al}}.
\end{multline}
Since $\be=\al/(\al-1)$, let note that $ h(x)^{\be/(\be-2)}=x^\al$ if $\ABS{x}\geq 1$. Then we fix $\eta=1/2$. And note
$$\A=\log{\int \exp\PAR{\frac{1}{2} h^{\be/(\be-2)}}d\mu_\al}<\infty.$$
Then we now fix $\e=\inf\BRA{\be/(\A(\be-2) 4{C_h} ),\be/((\be-2)4 {C_h})}$. We obtain
$$\ent{\mu_\al}{{(f-2)_+}^2}\leq\frac{{C_h}}{\e^{(\be-2)/2}} \int_{f\geq2}\ABS{\frac{f'}{f}}^\be f^2d\mu_\al
+\frac{1}{4}\ent{\mu_\al}{f^2}+\frac{1}{4}.$$
$\ent{\mu_\al}{f^2}\geq 1$, implies
\begin{equation*}
\ent{\mu_\al}{f^2}\leq 20\var{\mu_\al}{f}+
\frac{2{C_h}}{\e^{(\be-2)/2}} \int_{f\geq2}\PAR{\frac{f'}{f}}^\be f^2d\mu_\al.
\end{equation*}
\proofend

\begin{erem}
As we can note it in the demonstration:  the constant $A$ is universal
and the constant $B$ depends on the measure studied.  We will see in section~\ref{section-autre}
 that one can adapt the demonstration for other measures.
\end{erem}

\begin{erem}
\label{prop-grand}
With the same method as developed in
Proposition~\ref{prop-grand2} we can prove the inequality~\eqref{eq-grand} without
$\var{\mu_\al}{f}$.
Suppose that $\al\in]1,2]$.
There exists $A>0$ such that for any functions $f$ satisfying
$$\int f^2 d\mu_\al=1\,\, \text{and} \,\,\ent{\mu_\al}{f^2}\geq 1$$
 we have
$$
\ent{\mu_\al}{f^2}\leq
A\int\ABS{\frac{f'}{f}}^\be f^2d\mu_\al.
$$
\end{erem}

{\it Small entropy case}.

Let us now give the result when $\ent{\mu_\al}{f^2}$ is small.

\begin{eprop}
\label{prop-petit}
Let $\al\in[1,2]$.
There exists $A,A'>0$ such that for any functions $f\geq0$ satisfying
$$\int f^2 d\mu_\al=1\,\, \text{and} \,\,\ent{\mu_\al}{f^2}\leq 1$$
 we have
$$
\ent{\mu_\al}{f^2}\leq 34\var{\mu_\al}{f}+
A\int_{f\geq 2}\ABS{\frac{f'}{f}}^\be f^2d\mu_\al.
$$

when $\al\in]1,2]$, and if $\al=1$ we get
$$
\ent{\mu_\al}{f^2}\leq A'\var{\mu_\al}{f}.$$
\end{eprop}


{\noindent {\emph{\textbf{Proof of Proposition~\ref{prop-petit}}}}\\~}
\proofbegin~
Let $f\geq0$ satisfying $\int f^2d\mu_\al=1$.
Like in Proposition~\ref{prop-grand2}, we start with inequality (\ref{triv-ineq}), which readily implies
\begin{eqnarray}
\label{eq-dep}
\ent{\mu_\al}{f^2}=\int f^2\ln f^2d\mu_\al \leq 10\var{\mu_\al}{f}+ \int (f-2)_+^2\ln f^2d\mu_\al.
\end{eqnarray}

We will now control the second term of the right hand side of this last inequality via the use of
$\Phi$-Sobolev inequalities, namely Lemma \ref{lem-g}. Therefore we have to construct a function $g$,
greater than 2, which satisfies (for a well chosen $T$), when $\ent{\mu_\al}{f^2}\le1$,
\begin{description}
\item{(g1)} $\displaystyle \int_T^\infty g^2d\mu_\al\le 12$;
\item{(g2)} $\displaystyle\int_T^\infty (g-2)^2_+\Phi(g^2)d\mu_\al \ge C \int (f-2)_+^2\ln f^2d\mu_\al$;
\item{(g3)}$ \displaystyle\int_T^\infty g'^2d\mu_\al\le C\int_{[T,\infty[\cup
\{f\ge 2\}}\Psi\left(\left|{f'\over f}\right|\right)f^2d\mu_\al +D ~\ent{\mu_\al}{f^2}$,
\end{description}

with $\Phi(x)=\ln^{2(\al-1)\over\al}(x)$, $0<D\leq1/2$ and $\Psi(x)=x^\beta.$

\bigskip

Let now define $T_1<0$ and $T_2>0$ such that
$$
\mu_\al\PAR{]\infty,T_1]}=\frac{3}{8},\,\,\mu_\al\PAR{[T_1,T_2]}=\frac{1}{4}
\,\,\text{and}\,\,\mu_\al\PAR{[T_2,+\infty[}=\frac{3}{8}.
$$
Since $\int f^2d\mu_\al=1$ there exists $T\in [T_1,T_2]$ such that $f(T)\leq2$.

Let us define $g$ on $[T_1,\infty]$ as follow
\begin{equation*}
g=
2+(f-2)_+\ln^\gamma f^2  \text{ on } [T,\infty[,
\end{equation*}
where $\gamma=(2-\al)/(2\al)$.

 Function $g$ satisfies $g(T)=2$ and $g(x)\geq2$ for all $x\geq T$.
Let now compute $\int_T^\infty g^2d\mu_\al$. We have
\begin{eqnarray*}
\int_T^\infty g^2d\mu_\al & \leq & 2\int_{T_1}^\infty 4d\mu_\al+2\int_{T_1}^\infty (f-2)_+^2\ln^{2\gamma} f^2d\mu_\al\\
& \leq & 4+2\int_{[T_2,\infty[\cap\BRA{f\geq2}}f^2\ln^{2\gamma} f^2d\mu_\al.
\end{eqnarray*}
Since $2\gamma\in[0,1]$ we have $\ln^{2\gamma} f^2\leq\ln f^2$ on $\BRA{f\geq2}$. Then we obtain
by Lemma~\ref{lem-var}.iii
\begin{eqnarray*}
\int_T^\infty g^2d\mu_\al& \leq  &5+2\int_{{f\geq2}}f^2\ln^{2\gamma} f^2d\mu_\al\\
& \leq & 5+8\ent{\mu_\al}{f^2}\\
& \leq & 13,
\end{eqnarray*}
since $ \ent{\mu_\al}{f^2}\leq 1$.

Assumptions on Lemma~\ref{lem-g} are satisfied, we obtain by inequality~\eqref{eq-lemme}

\begin{equation*}
\int_{T}^\infty (g-2)_+^2\ln^{\frac{2(\al-1)}{\al}}g^2d\mu_\al
\leq {C_g}\int_{[T,\infty[\cap\BRA{g\geq2}} {g'}^2d\mu_\al.
\end{equation*}
Let us compare the various terms now.

Firstly let note $u={2(\al-1)}/{\al}$, we have
$$
(g-2)_+^2\ln ^u g^2=(f-2)_+^2\ln^{2\gamma}f^2\ln^u\PAR{2+(f-2)_+\ln^\gamma f^2}^2.
$$
On $\BRA{f\geq 2}$ we have
$2+(f-2)_+\ln^\gamma f^2\geq 2+(f-2)_+K,$ where $K=\ln^\gamma 4$.
Since $K\geq1$ and $u+2\gamma=1$ one has
$$
(g-2)_+^2\ln ^u g^2\geq(f-2)_+^2\ln^{2\gamma+u}f^2=(f-2)_+^2\ln f^2.
$$
Then we obtain
\begin{equation}
\label{eq-first}
 \int_{T}^\infty (f-2)_+^2\ln f^2d\mu_\al
\leq \int_{T}^\infty (g-2)_+^2\ln ^{\frac{2(\al-1)}{\al}}g^2d\mu_\al.
\end{equation}

Secondly one has on $\BRA{f\geq2}$
$$
g'=f'\ln^\gamma f^2\PAR{1+\gamma2^\gamma\frac{f-2}{f\ln f^2}},
$$
then using $\ln f^2\geq\ln 4$ one obtain
$$
\ABS{g'}^2\leq\ABS{f'}^2\ln^{2\gamma}f^2\PAR{1+\frac{\gamma2^\gamma}{\ln 4}}^2.
$$
Let note $D=\PAR{1+{\gamma2^\gamma}/{\ln 4}}^2$, one has
\begin{equation}
\label{eq-secondly}
\int_{[T,\infty[ \cap \BRA{f\geq 2}}g'^2d\mu_\al
\leq D\int_{[T,\infty[ \cap \BRA{f\geq 2}}f'^2\ln^{2\gamma}f^2d\mu_\al,
\end{equation}
on $[T,\infty[ \cap \BRA{f\geq 2}.$

Then, using inequalities~\eqref{eq-first} and~\eqref{eq-secondly}, there exists $C\geq0$
(independent of $T\in[T_1,T_2]$), such that
$$
\int_{T}^\infty (f-2)_+^2\ln f^2d\mu_\al\leq C\int_{[T,\infty[ \cap
\BRA{f\geq 2}} f'^2\ln^{2\gamma}f^2d\mu_\al.
$$

When $\al=1$, one has trivially that the right hand side of this last inequality is bounded, when $|f'|\le 1$,
 by $C\int_{f\ge2} f^2d\mu_\al$ which is itself bounded, by Lemma \ref{lem-var}.ii, by $8C\var{\mu_\al}{f}$ which
concludes the proof in this case.

When $\al\in]1,2]$, we apply Inequality~\eqref{eq-young} with $q=\al/(2-\al)$ and $p=\al/(2(\al-1))$.
We obtain  for every $\e>0$,
\begin{multline*}
\int_{[T,\infty[ \cap \BRA{f\geq 2}} \PAR{\frac{f'}{f}}^2\PAR{\ln^{2\gamma} f^2}f^2d\mu_\al\leq \\
\frac{2(\al-1)}{\al\e^{\frac{2-\al}{2(\al-1)}}}\int_{[T,\infty[ \cap \BRA{f\geq 2}}
\ABS{\frac{f'}{f}}^\be f^2d\mu_\al +
\e\frac{2-\al}{\al}\int_{[T,\infty[ \cap \BRA{f\geq 2}} f^2 \ln f^2d\mu_\al.
\end{multline*}
Fix $\e$ such that $\e C\frac{2-\al}{\al}<1/16$, then  there exists $A>0$ such that
$$
\int_{T}^\infty (f-2)_+^2\ln f^2d\mu_\al\leq
A\int_{[T,\infty[ \cap \BRA{f\geq 2}} \ABS{\frac{f'}{f}}^\be f^2d\mu_\al +
\frac{1}{16}\int_{[T,\infty[ \cap \BRA{f\geq 2}} f^2 \ln f^2d\mu_\al.
$$
Using Lemma~\ref{lem-var}.iii we have,
\begin{equation*}
\int_{T}^\infty (f-2)_+^2\ln f^2d\mu_\al\leq\\
A\int_{[T,\infty[ \cap \BRA{f\geq 2}} \ABS{\frac{f'}{f}}^\be f^2d\mu_\al +
\frac{1}{4}\ent{\mu_\al}{f^2}.
\end{equation*}

The same method can be used on $]-\infty,T]$ and then , there is $A'>0$ such that
\begin{equation*}
\int_{-\infty}^T (f-2)_+^2\ln f^2d\mu_\al\leq\\
A'\int_{[-\infty,T] \cap\BRA{f\geq 2}} \ABS{\frac{f'}{f}}^\be f^2d\mu_\al +
\frac{1}{4}\ent{\mu_\al}{f^2}.
\end{equation*}

And then, we get
$$
\int(f-2)_+^2\ln f^2d\mu_\al\leq\\
(A+A')\int_{{f\geq 2}} \PAR{\frac{f'}{f}}^\be f^2d\mu_\al +
\frac{1}{2}\ent{\mu_\al}{f^2}.
$$
Note that constants $A$ and $A'$ don't depend on $T$.

Then, by inequality~\eqref{eq-dep}, Proposition~\ref{prop-petit} is proved.
\proofend

\medskip

Let us give now a proof of the theorem.

{\noindent {\emph{\textbf{Proof of Theorem~\ref{thm-the}}}}\\~}
\proofbegin~
The proof of the theorem is a simple consequence of Propositions~\ref{prop-grand2} and  \ref{prop-petit}.
\proofend

\section{Extension to other measures}
\label{section-autre}

We will present in this section modified logarithmic Sobolev inequality of function $H$ for more general
measure than $\mu_\al$ which can be derived using the proof carried on in Section~\ref{sec-ex}: the large
entropy case where the optimal Hardy function $h$ is identified and used to derive the optimal $H$,
 and the small entropy case where $\Phi$ and $g$ (used on the proof of Proposition~\ref{prop-grand})
have to be identified leading to the same $H$ function.\par
\vspace{5pt}

Let us first consider the following probability measure
$\mu_{\al,\beta}$ for $\alpha\in[1,2]$ and $\beta\in\R$ defined by

$$\mu_{\al,\beta}(dx)=\frac{1}{Z}e^{-\phi(x)}dx\text{ where } \phi(x)=|x|^\al (\log|x|)^\beta \text{ for } |x|\ge 1$$

and $\phi$ twice continuously differentiable.

\begin{ethm}
There exists $A,B>0$ such that the measure $\mu_{\al,\beta}$ satisfies the
 following logarithmic Sobolev inequality: for any smooth $f$ on $\R$ such that $\int f^2d\mu_{\al\beta}=1$ and $f\geq0$,
 we have
\begin{equation}
\label{eq-hardy3}
\ent{\mu_{\al,\beta}}{f^2}\leq A\var{\mu_{\al,\beta}}{f}+
B\int_{{f}\geq 2}H\left(\ABS{\frac{f'}{f}}\right) f^2d\mu_{\al,\beta},
\end{equation}

where $H$ is positive smooth and given for $x\ge 2$ by
$$H(x)=\frac{x^{\alpha\over \alpha-1}}{\log^{\beta\over \alpha -1}x} \,\text{ if }\, \al\in]1,2[,\be\in\dR,$$
$$H(x)=x^2 e^{x^{1/\be}}\,\text{ if }\, \al=1,\be\in\dR^+ \,
\text{ and }\,
 H(x)=x^2 \log^{-\be}(x)\,\text{ if }\, \al=2,\be\in\dR^-.$$
\end{ethm}
\begin{eproof}
We will mimic closely the proof given in the $\mu_\al$ case, considering large and small entropy case. We will not
present all the calculus but give the essential arguments.\par\vspace{5pt}

Let now treat the case $\al\in]1,2[$.

{\it Large entropy.} We will first apply Lemma \ref{lem-h} to measure $\mu_{\al,\beta}$, one has then that
 $b_+$, $b_-$, $B_+$, $B_-$ are finite if one take $h$ positive smooth
$$
h(x)=\frac{x^{2-\al}}{\log^\beta x} \qquad |x|\ge 2.
$$
One has then to determine $H$ to construct $\psi$ such that there exists $\eta>0$ with $\eta\psi(h)$
exponentially integrable with respect to $\mu_{\al,\beta}$ and $H=\psi^*(x^2)$ where $\psi^*$ is the
Fenchel-Legendre transform of $\psi$.

Considering the exponential integrability condition leads us to consider $\psi(x)$ behaving asymptotically
as $x^{\al\over 2-\al}\log^{2\beta\over2-\al}x$. One may thus derive the asymptotic behavior of $\psi^*$ and
 finally $H$.\par\vspace{5pt}

 {\it Small entropy.} One desires here to apply Lemma
\ref{lem-g}, evaluating $\Phi$ and then build the
 function $g$ satisfying conditions (g1), (g2) and (g3). By Hardy's inequality and arguments in the proof
 of Lemma \ref{lem-g}, one may choose $\Phi$ for $x$ large enough as
$$\Phi(x)=\log^{2\frac{\al-1}{\al}}(x)\left(\log\log x\right)^{\frac{2\beta}{\al}}.$$
Setting then
$$g=2+(f-2)_+\log^{2-\al\over2\al}f^2(\log\log f^2)^{-\beta\over \al},$$
one may then verify (g1), (g2) and (g3) with $\Psi=H$ defined in the large entropy step.

Now if $\al=1$ and $\be\geq0$, then the same arguments gives that
for large enough $x$
$$\psi(x)=x\log^{2\be}x,\,\,\,\psi^*(x)=xe^{x^{1/(2\be)}}\,\text{
and }\,H(x)=x^2 e^{x^{1/\be}}.$$

 If  $\al=2$ and $\be\leq0$, we have for large enough $x$.
$$\psi(x)=\frac{1}{x}e^{2x^{-1/\be}},\,\,\,\psi^*(x)=x\log^{-\be} x\,\text{
and }\,H(x)=x^2 \log^{-\be} x.$$
\end{eproof}

\begin{erem}
\begin{enumerate}
\item Using once again Herbst's argument, we may derive concentration properties for
the measure $\mu_{\al,\beta}$ of desired order,  for every
Lipschitz function $F$ with $\| F\|_{lip}\le 1$, there exists $C>0$ such that, for all $\la>0$,
$$
\mu_{\al,\beta}\left(\ABS{F-\mu_{\al,\beta}(F)}\ge \lambda\right)
\le2 e^{-C\min(\lambda^\al\log^\beta\lambda,\la^2)}.
$$
\item Note that the Lata{\l}a-Oleszkiewicz inequalities $I(r)$ (see \cite{latala}) are not well
adapted for the family of measures $\mu_{\al,\beta}$. Indeed, using Hardy's characterization
of this inequalities obtained by Barthe-Roberto \cite[Th. 13 and Prop. 15]{barthe-roberto},
one may show that $\mu_{\al,\beta}$ satisfies an $I(\alpha/2)$ inequality if $\beta\ge0$ and an
 $I(\al/2-\epsilon)$ ($\epsilon$ being arbitrary small) for $\beta<0$, which entails consequently
not optimal concentration properties.
\item
By the characterization of the spectral gap property on $\R$, one obtains that each measure
 $\mu_{\al,\beta}$ satisfies a Poincar\'e inequality and thus a modified logarithmic Sobolev inequality.
\end{enumerate}
\end{erem}

Following the previous proof, we may generalize the family $\mu_{\al,\beta}$ adding an explicit multiplicative term
to the potential $|x|^\al \log^\beta |x|$, as for example $\log\log^\gamma |x|$ which will give us new modified
 logarithmic Sobolev inequality, but each of this new measure has to be considered ``one-by-one'' (we hope some general results
for $\phi$ convex). We may now
 state a result enabling us to get the stability of these modified logarithmic Sobolev inequality by addition
 of an unbounded perturbation: consider the measures
$$
d\tau_\al(x)=\exp\left(-|x|^\al-|x|^{\al-1}\cos(x)\right)\frac{dx}{Z_\al},\qquad \al\in]1,2],
$$
$$
d\gamma_{\al,b}(x)=(1+x)^be^{-x^\al}\frac{dx}{Z_{\al,b}}1_{x\geq0},\qquad \qquad \al\in]1,2],b\in\dR.
$$

\begin{eprop}
There exists $a>0$ such that the measures $\tau_\al$ and $\gamma_{\al,b}$ satisfy a logarithmic Sobolev
 inequality of function $H_{a,\al}$.
\end{eprop}

\begin{eproof}
Following the proof given in Section~\ref{sec-ex} , one sees that the result hold true once one may verify that the
Hardy's inequalities of Lemma \ref{lem-h} and Lemma \ref{lem-g} hold with the $h$ and $\Phi$ obtained for
the case of $\mu_\al$. It is easily checked once remarked that
$$\log \frac{d\tau_\al(x)}{dx}\sim_\infty -|x|^\al \qquad {\rm and}\qquad \left(\log \frac{d\tau_\al(x)}{dx}\right)'\sim_\infty
-(\al-1)|x|^{\al-1}$$
and the same for $\ga_{\al,b}$.
\end{eproof}



\newcommand{\etalchar}[1]{$^{#1}$}

\end{document}

%% file: macros.tex
\newcommand{\disp}{\displaystyle}

\newcommand{\dN}{\ensuremath{\mathbb{N}}}

\newcommand{\dR}{\ensuremath{\mathbb{R}}}

\newtheorem{ethm}{Theorem}[section]

\newtheorem{ecor}[ethm]{Corollary}

\newtheorem{eprop}[ethm]{Proposition}

\newtheorem{elem}[ethm]{Lemma}

\newtheorem{edefi}[ethm]{Definition}

\newtheorem{erem}[ethm]{Remark}

\newcommand{\proofend}{~$\rhd$}
\newcommand{\proofbegin}{~$\lhd$}

\newenvironment{eproof}
               {\noindent {\emph{\textbf{Proof}}}\\\proofbegin~}
               {\proofend\\}

\newcommand{\p}[4]{{#3}\!\left#1{#4}\right#2}

\newcommand{\ABS}[1]{\ensuremath{{\left| #1 \right|}}} 
\newcommand{\PAR}[1]{\ensuremath{{\left(#1\right)}}} 
\newcommand{\BRA}[1]{\ensuremath{{\left\{#1\right\}}}} 
\newcommand{\NRM}[1]{\ensuremath{{\left\Vert #1\right\Vert}}} 
\newcommand{\abs}[1]{\ensuremath{\lvert #1 \rvert}}

\renewcommand{\phi}{\varphi}

\renewcommand{\geq}{\geqslant}

\newcommand{\varf}[1]{\mathbf{Var}_{#1}}
\newcommand{\entf}[1]{\mathbf{Ent}_{#1}}

\newcommand{\ent}[2]{\p(){\entf{#1}}{#2}}

\newcommand{\var}[2]{\p(){\varf{#1}}{#2}}

\def\disp{\displaystyle}

\newcommand{\A}{\ensuremath{\mathcal A}}

\newcommand{\R}{\dR}

\newcommand{\al}{\alpha}

\newcommand{\be}{\beta}

\newcommand{\ga}{\gamma}

\newcommand{\ep}{\epsilon}

\newcommand{\la}{\lambda}

\newcommand{\e}{\varepsilon}

